\documentclass[11pt]{article}
 \usepackage{amssymb}
 \usepackage{epsfig}
 \usepackage{graphicx}
 \usepackage{pict2e}
 \newtheorem{Lemma}{Lemma}
 \newtheorem{Proposition}[Lemma]{Proposition}
 
 \newtheorem{Theorem}[Lemma]{Theorem}

 \newtheorem{Corollary}[Lemma]{Corollary}

 \newcommand{\Reals}{{\mathbb{R}}}
 
 \newcommand{\Ints}{{\mathbb{Z}}}
 \newcommand{\Ex}{{\mathbb{E}}}

 \newcommand{\TT}{\mbox{${\mathcal T}$}}

 \newcommand{\EE}{\mbox{${\mathcal E}$}}
 \newcommand{\VV}{\mbox{${\mathcal V}$}}
 \newcommand{\GG}{\mbox{${\mathcal G}$}}
 \newcommand{\LL}{\mbox{${\mathcal L}$}}

 \newcommand{\eps}{\varepsilon}

 \newcommand{\bx}{{\mathbf x}}
 \newcommand{\by}{{\mathbf y}}

 \newcommand{\sfrac}[2]{{\textstyle\frac{#1}{#2}}}
 \newcommand{\qed}{\ \ \rule{1ex}{1ex}}

 \newcommand{\len}{\, \mathrm{len}}

 \newcommand{\proof}{{\bf Proof.\ }}

\begin{document}
 \title{Which Connected Spatial Networks on Random Points have Linear Route-Lengths?}
 \author{David J. Aldous\thanks{Department of Statistics,
 367 Evans Hall \#\  3860,
 U.C. Berkeley CA 94720;  aldous@stat.berkeley.edu;
  www.stat.berkeley.edu/users/aldous.  Aldous's research supported by
 N.S.F Grant DMS-0704159. }
}

 \maketitle
 
  \begin{abstract}
In a model of a connected network on random points in the plane, one expects that the mean length  of the shortest route between vertices at distance $r$ apart should grow only as $O(r)$ as $r \to \infty$, but this is not always easy to verify.  We give a general sufficient condition for such linearity, in the setting of a Poisson point process.  In a  $L \times L$ square, define a subnetwork $\GG_L$ to have the edges which are present regardless of the configuration outside the square; 
the condition is that the largest component of $\GG_L$ should contain a proportion $1 - o(1)$ of the vertices, as $L \to \infty$.  The proof is by comparison with oriented percolation.  We show that the general result applies to the {\em relative neighborhood graph}, and establishing the linearity property for this network immediately implies it for a large family of {\em proximity graphs}.  
  \end{abstract}
 \vspace{0.1in}

 \vspace{0.4in}



 \vspace{0.4in}


 \newpage
 \section{Introduction}
 \label{sec-INT}
 The phrase {\em random networks} encompasses a very wide range of mathematical models and real-world motivations. 
 We will consider the {\em spatial network} setting where the vertices are points in two-dimensional space.  
 In ordinary language, ``network" often carries the implication of being connected, so it is ironic the the two most classical models for (infinite)
 random spatial networks are not connected.  These are 
 
 (i) [{\em Random geometric  graph} \cite{penrose-RGG}].  Start with a Poisson point process of vertices; put an edge between vertices $i,j$  if the Euclidean distance $d(i,j)$ is less than some prescribed constant $c$.  
 
 (ii) [{\em Bond percolation} \cite{gri99}].  Start with the square lattice; retain or delete edges at random.
 
 A third class of models {\em is} connected:
 
 (iii) [{\em Small world} models  \cite{newman-small}].  Start with the usual square grid of vertices and edges; add extra edges $(i,j)$ with probability $p(i,j)$ for prescribed 
$p(\cdot)$.

However this third class of  models is designed for settings where the relevant notion of path-length is ``number of edges in the path".  We will be concerned with the ``purely geometric" setting where the notion of {\em route length} is ``sum of Euclidean lengths of edges in the route", so there is a 
{\em network distance} between any two vertices, defined as the minimum (over all routes connecting the vertices) route length.
In contexts where  networks are designed to be connected and network distance is a topic of interest, using mathematical models which are not {\em a priori} connected is both conceptually unsatisfying and technically complicated, in that one is forced to consider quantities such as 
``mean network distance conditioned on the two vertices being in the unique infinite component".

This paper is a technical contribution to  a program (see  \cite{me-spatial-1} for a survey) studying connected networks on random points, where we view the  
 infinite Poisson process as proxy for a large finite unstructured set of points. 
The main result concerns the very general class of {\em CIDRPP} networks that we will now specify, after apologizing for the ugly name, and after mentioning a terminological convention we find helpful.  
Visualizing road networks, we write {\em city} and {\em road} and {\em route} for the objects in a network model we are studying; and write 
 {\em vertex (site)} and {\em edge} and {\em path} for mathematical graphical objects constructed in the course of the proofs.

The {\em PP} denotes Poisson process, precisely the rate 1 per unit area Poisson point process used as the model for city positions.  
The {\em DR} denotes {\em deterministic rule} for roads;  given the configuration of cities, whether or not a straight road links the two cities is determined by some deterministic rule, which need not be ``local" (it may depend on the entire configuration) but is required only to be invariant ({\em I}) under translation. (Note that using simple {\em random} rules would tend to leave some cities isolated and hence disconnect the network).  The resulting IDRPP network is a translation- invariant process in the plane, and the 
 ``random geometric  graph" in (i) is perhaps the simplest example.   Finally we want the network to be connected ({\em C}), and adding that requirement gives the class of  CIDRPP networks. 
The best known example is perhaps the Delaunay triangulation \cite{MR1686530} on the PP; another example 
is the minimal spanning tree (MST) on the PP (see section \ref{sec-MST}).   But more interesting for our purposes is the 
{\em relative neighborhood graph} ({\em RNG}) \cite{jaromczyk} defined by 
\begin{quote}
there is a road between two cities $x,y$ if and only if there is no other city $z$ with 
$\max(d(z,x), d(z,y))  < d(x,y)$
\end{quote}
where $d$ denotes Euclidean distance.  
See Figure 6 for an illustration.  This particular network is interesting because (loosely speaking) it is the sparsest connected graph that can be defined by a simple local rule; more precisely (see section \ref{sec-discuss}) there is a family of {\em proximity graphs} which are supergraphs of the RNG, and therefore can only have shorter network distances.

The examples above are planar graphs, for which the notion of {\em route} is
unambiguous.  Our results will apply to non-planar networks, if we make the
convention that if a road $(x_1,x_2)$ crosses a road $(x_3,x_4)$  then there is a
route $(x_1,x_4)$ via the junction.  In most natural examples the DR is also {\em rotation}-invariant, but we do not need to assume that.

\subsection{The linearity property}
Now consider some CIDRPP network.  
How should we formalize the ``linearity" property
\begin{quote}
mean network distance  between vertices at Euclidean distance $r$ apart should grow only as $O(r)$ as $r \to \infty$.
\end{quote}
By familiar properties of the PP, conditional on cities at prescribed points $z^1, z^2$ in $\Reals^2$ with $z^1 -  z^2 = z$, the network distance between them is distributed as if we ``planted" cities at the origin and at $z$, other cities being distributed as the PP, and considered the network distance $T_z$ between the planted cities.  Here $T_z$ is random with 
$|z| \leq T_z < \infty$.   
So we could formalize the linearity property as
\begin{equation}
\limsup \frac{\Ex T_z}{|z|} < \infty  \mbox{ as } |z| \to \infty .
\label{lin-prop-0}
\end{equation}
Unfortunately it seems technically hard to deal with this definition at the level of generality we seek, because the deterministic rule
 for  edge presence may depend on the entire configuration.  Instead we shall use the ``integrated out" analog.  
 Denote by $\xi$ the cities in the PP, and
 write $\ell(\xi,\xi^\prime)$ for network distance.  
 Say a CIDRPP network has the {\em linearity property} if, for any choice of bounded subsets $A, B$ of $\Reals^2$,
\begin{equation}
\Ex \sum_{\xi \in A} \sum_{\xi^\prime \in z+B} \ell(\xi,\xi^\prime) = O(|z|) \mbox{ as } |z| \to \infty .
\label{lin-prop}
\end{equation}
Here $z + B :=\{z+z^\prime: \ z^\prime \in B\}$.  

Clearly it is enough to prove (\ref{lin-prop}) when $A$
and $B$ are  squares (of some one arbitrary given size) centered at the origin;  
by considering very small squares 
 we see  this definition is intuitively similar to (\ref{lin-prop-0}).

A glance at Figure 6 suggests that the RNG has this linearity property, and Monte Carlo simulations \cite{me-spatial-1} suggest a limit of about $1.38$ in (\ref{lin-prop-0}).
 But we do not know any proof which is ``constructive", in the sense of first describing an algorithm for choosing a route and then bounding the mean length of the chosen route.   Rather than seeking some non-constructive proof tied to this particular model it seemed more useful to seek a more widely applicable result, 
 and this consideration was the motivation for this paper.  
 It is clear (intuitively, at last -- see section \ref{sec-MST}) that the MST does not possess the linearity property, so we need to impose some further assumptions on a  CIDRPP network to ensure the linearity property.

It turns out that one condition is sufficient.
Consider the $L \times L$ square $[0,L]^2$. 
Then consider the subnetwork $\GG_L$ defined in words as
\begin{eqnarray}
\mbox{the cities in $[0,L]^2$, with the roads that are present} \nonumber\\
\mbox{ regardless of the
configuration of cities outside $[0,L]^2$.} \label{def-GG}
\end{eqnarray}
To say this more formally, consider a deterministic configuration of cities in the
plane, write $\bx$ for the configuration of cities inside $[0,L]^2$, and 
$x_1,x_2$ for two such cities, and write $\by$ for the configuration of cities
outside $[0,L]^2$.  The deterministic rule defining the network is an
$\{{\tt edge, no-edge}\}$-valued function 
$f(x_1,x_2, \bx, \by)$ specifying whether to put an edge between $x_1$ and $x_2$.  
Define a new function 
\[\bar{f}(x_1,x_2, \bx)  = {\tt edge} \mbox{ iff } f(x_1,x_2, \bx, \by) = {\tt edge}
\mbox{ for all } \by . \]
Now $\GG_L$ is the network obtained by applying rule $\bar{f}$ to assign roads to a
PP of cities inside $[0,L]^2$.  
By construction, $\GG_L$ is a subnetwork of the original network in the infinite plane.

The subnetwork $\GG_L$ need not be connected, so write $N^0_L$ for the number of
cities inside $[0,L]^2$ that are \underline{not} in the largest component of $\GG_L$.
Consider the {\em asymptotic essential connectedness} property
\begin{equation}
L^{-2} \ \Ex N^0_L \to 0 \mbox{ as } L \to \infty.
\label{AC}
\end{equation}
\begin{Theorem}
\label{T1}
If a CIDRPP network satisfies  the asymptotic essential connectedness property
(\ref{AC}) then it has the linearity property (\ref{lin-prop}).
\end{Theorem}
We prove this in section \ref{sec-PT1}  by comparison with oriented percolation.
This general method is well known in the theoretical spatial random processes literature,  and is indeed a central focus of the 1988 monograph of Durrett \cite{dur88}.

Proposition \ref{P2} in section \ref{sec-RNG} checks that the relative neighborhood graph satisfies the assumptions, and hence the conclusion, of Theorem \ref{T1}.

\subsection{Remarks}
Let us conclude the introduction with some brief remarks.

{\bf 1.} We should acknowledge that Theorem \ref{T1} is ``purely theoretical" in that the argument gives  very large bounds.  
Obtaining numerically reasonable bounds in networks like the RNG seems a hard problem.

{\bf 2.}  Our linearity property (\ref{lin-prop}) can be viewed as the average-case analog of the worst-case concept of a 
{\em spanner}, described in section \ref{sec-discuss}.

{\bf 3.} 
A counter-intuitive result \cite{me116} says that, if one seeks to design a network over a PP with only the two goals of minimizing 
$\Ex T_z$ (defined above (\ref{lin-prop-0})) for large $|z|$, and minimizing mean road length per unit area ($\alpha$, say), then one can construct networks such that\\
(i)  $|z|^{-1} \Ex T_z \to 1 \mbox{ as } |z| \to \infty$ \\
(ii) $\alpha$ is arbitrarily close to $\alpha_{\mbox{{\tiny ST}}}$, the value for the Steiner tree, the smallest value for any connected network.\\
Consequently,  while Theorem \ref{T1} is intended to be useful for establishing linearity in networks that are given to us, it isn't needed for networks that we are allowed to design.

{\bf 4.} Where (\ref{lin-prop-0}) holds it is natural to conjecture that for fixed $\theta$, taking 
$z = (r,\theta)$ in radial coordinates, the limit 
\[ 
\lim_{r \to \infty} \frac{\Ex T_{(r,\theta)} }{r} 
\]
always exists, and that this property 
(analogous to the {\em shape theorem} in first passage percolation)
should be easy to prove by subadditivity.  
But counter-intuitively, existence of the limit does not seem easy to prove; a result 
of this type, formalized by an analog of (\ref{lin-prop}),  will be given in a companion paper \cite{me-spatial-5}.

{\bf 5.} For technical reasons the hypothesis of \cite{me-spatial-5} uses a ``$L^2$" analog of the ``$L^1$" 
linearity property (\ref{lin-prop}), which in fact (section \ref{sec-completing}) follows from the proof of Theorem \ref{T1}: we record the general 
$k$'th moment form as 
\begin{Corollary}
\label{C1}
Under the hypotheses of Theorem \ref{T1}, for each $k \ge 1$ and each choice of 
bounded subsets $A, B$ of $\Reals^2$,
\begin{equation}
\sup_z  \frac{\Ex \sum_{\xi \in A} \sum_{\xi^\prime \in z+B} \ell^k(\xi,\xi^\prime)}
{\max(1,|z|^k) } < \infty  .
\label{lin-prop-k}
\end{equation}
\end{Corollary}

\section{Proof of Theorem \ref{T1}}
\label{sec-PT1}
The argument uses a very standard methodology in percolation theory: 
a  block construction and a  comparison with oriented percolation.
The closest explicit results in the percolation literature are those such as
 \cite{MR1404543} 
studying  network distance within the infinite component and giving bounds which hold for all $p$ in 
the supercritical regime.  
But several difficulties arise in applying such results directly in our setting.
Because we are bounding an expectation we can't simply ignore events of  probability $\to 0$; 
and we need information on distance-related quantities for all pairs of sites, not just the ``open" sites.

Anyway, instead of trying to exploit sophisticated results from percolation theory 
we shall give an argument using only one 
fundamental result,
Proposition \ref{P1} below.  
We outline the overall argument in section \ref{sec-outline}, after giving the block construction in 
section \ref{sec-block}.

\subsection{Comparison with oriented percolation}
\label{sec-OP}
Following Durrett \cite{MR757768} we set up a tilted square grid, the graph $\LL = (\VV,\EE)$ with vertices 
(sites) $\VV = \{ (n_1,n_2): n_i \in \Ints, n_1+n_2 \mbox{ even} \} $
and edges $(n_1,n_2) \ - \  (n_1 \pm 1, n_2+1)$. 
A {\em doubly-infinite oriented-up path}  is a path 
$\{(m(n),n), \ -\infty < n < \infty\}$ where $m(n+1) = m(n) \pm 1$. 

Fix $0<p<1$.  
A  {\em 1-dependent $p$-site process} is a random process in which each site $v \in
\VV$ is declared {\tt good} or {\tt bad} is such a way that 

\noindent
(i) $P(v\ \mbox{ is }  {\tt good}) =
p$

\noindent
(ii) for disjoint subsets $V_1, V_2$ of $\VV$ such that there is no edge $(v_1,v_2)$ with
$v_i \in V_i$,
the family of events $\{v \  \mbox{ is } {\tt good}\}_{v \in V_1}$ 
is independent of the family $\{v \ \mbox{ is }  {\tt good}\}_{v \in V_2}$.

\medskip
In such a process, a {\em good path} is a path through only good sites.  
Our argument is based on  the following standard result, illustrated in Figure 1.
\begin{Proposition}
\label{P1}
There exist $p_0<1$ and finite constants $(\kappa_k, k \geq 1)$ such that a 1-dependent $p$-site process with $p \geq p_0$ has
the following property.
Let $T \geq 2$ be the smallest even positive  integer such that site $(T,0)$  is in a
doubly-infinite oriented-up good path.  Then $\Ex T^k \leq \kappa_k$  for all $k \ge 1$.
\end{Proposition}
This result is given in Durrett  \cite{MR757768} section 10 by a ``contour argument". 
(In fact \cite{MR757768}  does the semi-infinite case, but the doubly-infinite case is similar.)

\setlength{\unitlength}{0.29in}
\begin{picture}(8,7)(-4,-1)
\put(0,0){\circle{0.53}}
\put(2,0){\circle*{0.53}}
\put(4,0){\circle*{0.53}}
\put(6,0){\circle*{0.53}}
\put(8,0){\circle*{0.53}}
\put(1,1){\circle*{0.53}}
\put(3,1){\circle*{0.53}}
\put(5,1){\circle*{0.53}}
\put(7,1){\circle*{0.53}}
\put(0,2){\circle*{0.53}}
\put(2,2){\circle{0.53}}
\put(4,2){\circle*{0.53}}
\put(6,2){\circle{0.53}}
\put(8,2){\circle*{0.53}}
\put(1,3){\circle*{0.53}}
\put(3,3){\circle*{0.53}}
\put(5,3){\circle{0.53}}
\put(7,3){\circle*{0.53}}
\put(0,4){\circle*{0.53}}
\put(2,4){\circle*{0.53}}
\put(4,4){\circle*{0.53}}
\put(6,4){\circle*{0.53}}
\put(8,4){\circle*{0.53}}
\put(1,5){\circle{0.53}}
\put(3,5){\circle{0.53}}
\put(5,5){\circle{0.53}}
\put(7,5){\circle*{0.53}}
\put(7.7,0.3){\vector(-1,1){0.4}}
\put(7.3,1.3){\vector(1,1){0.4}}
\put(7.7,2.3){\vector(-1,1){0.4}}
\put(6.7,3.3){\vector(-1,1){0.4}}
\put(6.3,4.3){\vector(1,1){0.4}}
\put(1.4,1.3){$(0,0)$}
\put(7.8,1.3){$(T,0)$}
\end{picture}

{\bf Figure 1.}  {\small Black discs denote good sites.  The arrows indicate  edges  of one doubly-infinite oriented-up path starting through $(T,0)$.}

\subsection{The block construction}
\label{sec-block}
Fix large $L$.
Write $\LL_L = (\VV_L,\EE_L)$ for the tilted square grid scaled by $L$, so its vertex-set is 
$\VV_L = \{ (Ln_1,Ln_2): n_i \in \Ints, n_1+n_2 \mbox{ even} \} $.
Write $\VV_L^*$ for the set of midpoints of edges of $\EE_L$.  
For each $v \in  \VV_L \cup \VV_L^*$  define the following objects.

$S_v$ is the open $L \times L$ square centered at $v$;

$\GG_v$ is the subnetwork (\ref{def-GG}) in $S_v$;  that is, 
 the roads within $S_v$ that are present regardless of the
configuration of cities outside $S_v$;

$N_1(\GG_v)$ is the number of cities in the largest component of $\GG_v$;

$\len(\GG_v)$ is the total length of the roads of $\GG_v$.

\medskip
\noindent
Note that the squares $(S_v)_{v \in \VV_L}$ are disjoint, 
as are the squares $(S_v)_{v \in \VV^*_L}$.
Call $v \in \VV_L$ and $v^* \in \VV_L^*$ {\em diagonally adjacent} if $v^*$ is the midpoint of some edge at $v$, 
in which case 
$S_v \cap S_{v^*}$ is  a $L/2 \times L/2$ subsquare of $S_v$.
See Figure 2.

\setlength{\unitlength}{0.05in}
\begin{picture}(50,50)(-15,-4)
\put(0,0){\circle*{2}}
\put(40,0){\circle*{2}}
\put(0,40){\circle*{2}}
\put(40,40){\circle*{2}}
\put(20,20){\circle*{2}}

\put(8.6,9.2){$\bigotimes$}
\put(28.6,9.2){$\bigotimes$}
\put(8.6,29.2){$\bigotimes$}
\put(28.6,29.2){$\bigotimes$}

\put(10,10){\line(1,0){20}}
\put(20,20){\line(1,0){20}}
\put(10,30){\line(1,0){20}}
\put(20,40){\line(1,0){20}}

\put(10,10){\line(0,1){20}}
\put(20,20){\line(0,1){20}}
\put(30,10){\line(0,1){20}}
\put(40,20){\line(0,1){20}}

\put(20,17){$v$}
\put(30,32){$v^*$}
\put(6.6,17){$S_v$}
\put(31,41){$S_{v^*}$}

\put(43,10){\vector(0,1){10}}
\put(43,10){\vector(0,-1){10}}
\put(45,10){$L$}

\put(-4,40){$w$}

\end{picture}

{\bf Figure 2.} 
{\small The $\bullet$ are vertices of $\VV_L$ and the $\bigotimes$ are vertices of $\VV^*_L$. 
Here $v$ and $v^*$ are diagonally adjacent.
}

\vspace{0.14in}
\noindent
Call $v \in  \VV_L \cup \VV_L^*$ {\em nice} if \\
(i) each of the four natural $L/2 \times L/2$ subsquares of $S_v$ has between $0.24L^2$ and $0.26 L^2$ cities ; \\
(ii) $N_1(\GG_v) \geq 0.99 L^2$; \\
(iii) $\len(\GG_v) \leq c_L $\\
where the constant $c_L$ will be specified later.  
Now call $v \in \VV_L$ {\em good} if $v$ and its four diagonally adjacent midpoints $v^* \in \VV^*_L$ are all nice.  

Using hypothesis (\ref{AC}), we can choose $L$ and $c_L$ such that $P(v\  \mbox{ is } {\tt good}) \geq
p_0$, where $p_0$ is the value in Proposition \ref{P1}.   
Let us emphasize that $L$ and $c_L$ remain fixed for the rest of the argument; we are never going to say ``let $L \to \infty$".

The construction above ensures the 1-dependent property for the process $\{v \  \mbox{ is } {\tt good}\}_{v \in \VV_L}$, because if 
$v_1$ and $v_3$ are not adjacent in $\EE_L$ then squares $S_{v^*_1}$ and $S_{v^*_3}$ (centered at midpoints 
$v^*_i$ of edges at $v_i$) cannot overlap.  

Observe the following consequence of parts (i) and (ii) of the definition of {\em nice} above. 
If $v^* \in \VV_L^*$ is a midpoint of an edge at $v \in \VV_L$ and if both $v$ and $v^*$ are nice, then the largest components of
$\GG_v$ and of $\GG_{v^*}$ must have a common city in the subsquare where the two squares intersect.
So for any cities $\xi$ and $\xi^*$ in the largest components of $\GG_v$ and of $\GG_{v^*}$, there exists a route in the underlying network from $\xi$ to $\xi^*$ for which 
 each road used is in $\GG_v$ or $\GG_{v^*}$.  
 (Note this is where the ``regardless of \ldots" in the definition of $\GG$ comes into play; this route exists, under the stated assumptions, regardless of other aspects of the configuration of all cities).
 
 Taking into account the definition of {\em good} site, we deduce the following.
Consider a good path $v_0,v_1,\ldots,v_q$ in $\VV_L$, and write $(v^*_i)$ for the midpoints of edges along that path. 
Then, for any cities $\xi_0$ and $\xi_q$ in the largest components of $\GG_{v_0}$ and $\GG_{v_q}$, there exists a route from $\xi_0$ to $\xi_q$ 
in which  each road used is in some $\GG_{v_i}$ or $\GG_{v_i^*}$.
Such a route has total length bounded by $(2q+1) c_L$, by (iii).  
Now for an oriented-up path we must have $q \le d(v_0,v_q)/L$;  also $d(v_0,v_q) \le d(x_0,x_q) + L \sqrt{2}$.
Recalling our notation $\ell(\xi,\xi^\prime)$ for network distance, we have shown
\begin{Lemma}
\label{Lxixi}
Let $\xi,\xi^\prime$ be cities in the PP, in the $L \times L$  squares centered at vertices $v, v^\prime$ of $\VV_L$, 
and in the largest components of $\GG_v$ and of $\GG_{v^\prime}$.  
If there exists an oriented-up good path in $\VV_L$ from $v$ to $v^\prime$ then 
\[ \ell(\xi,\xi^\prime) \le (4 + 2 L^{-1}d(\xi,\xi^\prime)) c_L . \]
\end{Lemma}

\subsection{Outline of argument}
\label{sec-outline}
In outline, the argument for Theorem \ref{T1} is rather clear.  \\
(i)  Given a small square $A$ at the origin and the translated square $z + A$ for large $|z|$,  
Lemma \ref{Lxixi} and Proposition \ref{P1} allow us to find two routes, one through some city $\xi$ near $A$ and the other through some city $\xi^\prime$ near $z +A$,  which can be oriented so as to meet
somewhere in between, and so create a route from $\xi$ to $\xi^\prime$ of length $O(|z|)$. \\
(ii) We want to make a route from cities in  $A$ to $\xi$, with  route-length $O(1)$. 
We do this by constructing in $\VV_L$ a circuit, with length $O(1)$,  of good sites including the site $v$ near $\xi$ and 
encircling the origin.  
This implies existence of a length $O(1)$ encircling circuit through $\xi$ in the underlying network.  Connectivity now ensures a length $O(1)$ route from cities in  $A$ to $\xi$.

Saying this carefully with the specific ``tilted square grid" construction of section \ref{sec-block} 
leads to technical difficulties in making paths meet.
We deal with this difficulty by  creating paths which can be oriented in more narrowly specified directions, by using a number of
``affine lattices" described next.  
Then part (i) is formalized in section \ref{sec-global} and part (ii) in section \ref{sec-local}.  
Where the section \ref{sec-block} arguments carry over straightforwardly to the ``affine lattice" setting we will omit the details.

\newpage
\subsection{Affine lattices}

\setlength{\unitlength}{0.3in}
\begin{picture}(10,10)(-3,-5.4)
\put(0,0){\circle*{0.2}}
\put(0,3){\circle*{0.2}}
\put(1,2){\circle*{0.2}}
\put(2,1){\circle*{0.2}}
\put(3,0){\circle*{0.2}}
\put(2,-1){\circle*{0.2}}
\put(1,-2){\circle*{0.2}}
\put(0,-3){\circle*{0.2}}

\put(0.2,2.87){(0,3)}
\put(1.2,1.87){(1,2)}
\put(2.2,0.87){(2,1)}
\put(0.2,-3.11){(0,-3)}
\put(3.2,-0.11){(3,0)}

\put(0,3){\line(0,-1){6}}
\put(0,0){\line(1,2){1}}
\put(0,0){\line(2,1){2}}
\put(0,0){\line(1,0){3}}
\put(0,0){\line(1,-2){1}}
\put(0,0){\line(2,-1){2}}

\put(7,0){\circle*{0.2}}
\put(7,3){\circle*{0.2}}
\put(8,2){\circle*{0.2}}
\put(9,1){\circle*{0.2}}
\put(10,0){\circle*{0.2}}
\put(8,-1){\circle*{0.2}}
\put(9,-2){\circle*{0.2}}
\put(10,-3){\circle*{0.2}}

\put(7,0){\line(0,1){3}}
\put(7,0){\line(1,2){1}}
\put(8,-1){\line(0,1){3}}
\put(8,-1){\line(1,2){1}}
\put(9,-2){\line(0,1){3}}
\put(9,-2){\line(1,2){1}}
\put(10,-3){\line(0,1){3}}

\put(7,-3){\circle*{0.2}}
\put(8,-4){\circle*{0.2}}
\put(9,-5){\circle*{0.2}}

\put(7,-3){\line(0,1){3}}
\put(7,-3){\line(1,2){1}}
\put(8,-4){\line(0,1){3}}
\put(8,-4){\line(1,2){1}}
\put(9,-5){\line(0,1){3}}
\put(9,-5){\line(1,2){1}}

\put(0.3,1.6){a}
\put(7.3,1.6){a}
\put(0.9,1.1){b}
\put(1.8,0.4){c}
\put(1.8,-0.6){d}
\put(0.9,-1.5){e}
\put(0.3,-1.9){f}
\end{picture}

{\bf Figure 3.}   
{\small The six sectors (left) and the affine lattice $\LL^a = (\VV^a,\EE^a)$ associated with sector $a$ (right).
}

\vspace{0.12in}
\noindent
Divide the half-plane into six sectors, labeled $a - f$, as on the left of Figure 3.
Associate with each sector an ``affine lattice", written as e.g. 
$\LL^a = (\VV^a,\EE^a)$, where the case ``$a$" is illustrated on the right of Figure 3.
For $a,b,c$ the vertex-set is 
 $\VV = \{ (n_1,n_2): n_i \in \Ints, n_1+n_2 \mbox{ divisible by } 3\} $.
The edge-set illustrated for $\EE^a$ has edges from $(n_1,n_2)$ to 
$(n_1+ 0, n_2+3)$ and to $(n_1 + 1, n_2 +3)$, 
where the added vectors $(0,3)$ and $(1,2)$ are the points (Figure 3, left) 
defining the boundary of sector $a$; analogously for $\LL^b$ and $\LL^c$.
The affine lattices $\LL^f,\LL^e,\LL^d$ can be viewed as the 
top-bottom reflections of 
$\LL^a,\LL^b,\LL^c$, and their vertex set is
$ \{ (n_1,n_2): n_i \in \Ints, n_1- n_2 \mbox{ divisible by } 3\} $.

We will describe a construction  involving $\LL^a$; the same construction can be applied to each  of these lattices.

We can obtain $\LL^a$ via a linear transformation applied to the tilted square grid 
$\LL = (\VV,\EE)$ of section \ref{sec-OP}.
As at the start of section \ref{sec-block}, scale $\LL^a$ by $L$ to get a lattice 
$\LL^a_L = (\VV^a_L,\EE^a_L)$.

We now repeat the construction in section \ref{sec-block},
replacing squares by their images (parallelograms) under the linear transformation.
This involves some obvious changes to constants, which we will not write out in
detail:  for instance the ``$L^2$" in the definition  of {\em nice} is replaced by 
$L^2$ times the area of the appropriate basic parallelogram.
Making the analogous definition of
$v \in \VV^a_L$ being {\em good} for $\LL^a_L$,
we can copy the argument for Lemma \ref{Lxixi} to obtain 
\begin{Lemma}
\label{Lpaths}
There exists a constant $A$ such that the following property, stated for $\LL^a_L$, holds for each of the six affine lattices.
Let $\xi,\xi^\prime$ be cities in the PP, in the scaled parallelograms 
$S_v^a$ and $S^a_{v^\prime}$ centered at vertices $v$ and $v^\prime$ of $\VV^a_L$, 
and in the largest components of the corresponding $\GG^a_v$ and  $\GG^a_{v^\prime}$.  
If there exists an oriented good path in $\VV^a_L$ from $v$ to $v^\prime$ then 
\[ \ell(\xi,\xi^\prime) \le A (1 + d(\xi,\xi^\prime))  . \]
\end{Lemma}
{\em Remarks on Lemma \ref{Lpaths}.}
(i) In the bound we are now supressing the dependence on $L$, which is fixed. \\
(ii) A path in (say) $\LL^a_L$ is {\em oriented-right} if each edge is 
a translate of one of the vectors bounding the sector in Figure 3, and 
{\em oriented-left} if it is a reversal of such, and {\em oriented} if either 
 oriented-right or oriented-left. \\
(iii) Hypothesis (\ref{AC}), the asymptotic essential connectedness property,
was stated for unit squares scaled by $L$, whereas to prove Lemma \ref{Lpaths} we need it to hold
for a fixed parallelogram $\diamondsuit$ scaled by $L$.  But this is true because for
any $\eps$  we can find a finite collection of (overlapping) squares inside $\diamondsuit$ whose 
union has area at least $(1 - \eps) \times $ the area of $\diamondsuit$. \\

Proposition \ref{P1} transfers directly to the affine lattice setting.
To state the conclusion, recall we are interested in the distance from site 
$(0,0)$ to a nearby site in (say) $\LL^a_L$ which is in a doubly-infinite oriented
path. 
We measure distance by considering the nearest such site in a specified direction
along 
either the reverse diagonal 
$\{(n_1, n_2): n_2 = - n_1\}$ 
or the forward diagonal 
$\{(n_1, n_2): n_2 =  n_1\}$, 
choosing a diagonal which does not bisect the sector under consideration.
Write the Euclidean distance from $(0,0)$ to the nearest distinct such site as 
e.g. $T^a_\nwarrow$, the arrow indicating the diagonal and direction in which 
we are searching.  
Call the site 
$v^a_\nwarrow$.
There is a finite collection $\TT$ of random variables such as $T^a_\nwarrow$ 
arising from different choices of sector and diagonal direction, and Proposition \ref{P1} implies
there exist
finite constants $\kappa_k^\prime$ such that for each $T \in \TT$,
\begin{equation}
\Ex T^k \le \kappa_k^\prime , \quad 1 \le k < \infty .
\label{Tkbound}
\end{equation}

\subsection{Global routes}
\label{sec-global}
The construction is illustrated in Figure 4. 
Consider the origin $(0,0)$ and some site $v$ in the scaled square lattice $\Ints^2_L$.
Then $v$ is in some sector, say sector $b$.   
Consider the two adjacent sectors, in this case $a$ and $c$.  
In the figure, 
$v_0 = v^a_\nwarrow$ is the closest site to $(0,0)$ along the diagonal $\nwarrow$ which is in a 
doubly-infinite oriented path in $\LL^a_L$, 
and 
$v_1$ is the closest site to $v$ along the diagonal $\nwarrow$ which is in a 
doubly-infinite oriented path in $\LL^c_L$.
One can always choose the diagonal and direction so that, as in Figure 4, 
for each edge $e_0$ in the path from $v_0$ toward the crossing point $z$ 
 and each edge $e_1$ in the path from $v_1$ toward the crossing point $z$, the edges are at angles bounded below 
$\pi/2$.
So the length of the path from $v_0$ to $v_1$ via the crossing point, 
interpreting the lattice edges as line segments in $\Reals^2$, is 
$O( |v| + |v_0| + d(v,v_1))$, 
by elementary geometry considerations.  
Then as in the argument for Lemmas \ref{Lxixi} and \ref{Lpaths},
we can associate with this lattice path a route in the underlying network,
of length at most $c_L \times$ lattice path length.
This construction leads to the following bound, discussed below.

\setlength{\unitlength}{0.09in}
\begin{picture}(50,40)(-20,-7)
\multiput(0,0)(-1,1){10}{\circle*{0.4}}
\multiput(-4,-2)(1,5){6}{\line(0,1){3}}
\multiput(-5,-4)(1,5){6}{\line(1,2){1}}

\multiput(21,19)(-1,1){10}{\circle*{0.4}}
\multiput(29,23)(-5,-1){9}{\line(-3,0){3}}
\multiput(26,23)(-5,-1){9}{\line(-2,-1){2}}

\put(0.3,-0.4){(0,0)}
\put(21.5,18.7){$v$}
\put(19.5,20.7){$v_1$}

\put(-1.5,17.5){$z$}

\put(-5,2.4){$v_0$}

\end{picture}

{\bf Figure 4.}
{\small  The ``global" part of route construction.
}

\vspace{0.12in}
\noindent

\begin{Proposition}
\label{P10}
Given any $v \in \Ints^2_L$, there exist sites $v_0, v_1 \in \Ints_L^2$ 
and sectors, say $a$ and $c$, such that 
$v_0$ (resp. $v_1$) is a good site for $\LL^a_L$ (resp. $\LL^c_L$), 
and for any cities 
$\xi,\xi^\prime$ in the scaled parallelograms 
$S_{v_0}^a$ and $S^c_{v_1}$ centered at vertices $v_0$ of $\VV^a_L$ and $v_1$ of $\VV^c_L$,
and in the largest components of the corresponding $\GG^a_{v_0}$ and  $\GG^c_{v_1}$,  
\[ \ell(\xi,\xi^\prime) \leq B(1 + d((0,0),v_0) + d(v, v_1) + |v|) . \]  
Here $B$ is a constant and 
\begin{equation}
\mbox{ each of $d((0,0),v_0)$ and $d(v, v_1)$ is distributed as some $T \in \TT$.}
\label{dvdv}
\end{equation}
\end{Proposition}

{\em Discussion.}
One notational detail is that the affine lattices $\LL^c_L$ only contain 
as vertices $1/3$ of the vertices of $\Ints^2_L$, so $v$ may be 
in one of the two translated copies of $\LL^c_L$ which cover the other vertices; 
this does not affect the argument.

An issue of more substance is to check that, where the two lattice paths cross, we
can connect one associated route in the underlying network to the other route.  
To check this, consider Figure 2. 
Any point $z$ within an edge of the lattice is between some $v \in \VV_L$ and a
diagonally adjacent
$v^* \in \VV^*_L$. 
Suppose, as in Figure 2, that $v^*$ is closer than $v$ to $z$. 
Properties (i) and (ii) of {\em nice} imply that 
$S_{v^*}$ contains at most $1.04L^2$ cities, and that the largest component of 
$\GG_{v^*}$ contains at least $0.99 L^2$ cities.
Transforming to the affine lattice setting, there are numerical constants 
$\alpha_i, \beta_i = (0.99/1.04) \alpha_i, \ i = a,\ldots, f$, and
 the crossing point $z$ in Figure 4 is in some 
parallelogram $S^a$ for which \\
(i)   $S^a$ contains at most $\alpha_a L^2$ cities; \\
(ii) the largest component of the associated $\GG^a$ contains at least $\beta_a L^2$ 
cities;\\
(iii) $z$ is on the line from the center to a corner of $S^a$, closer to the center than to the corner;\\
and similarly for another parallelogram $S^c$.
These constraints force the intersection of the largest 
components of $\GG^a$ and of $\GG^c$ to be non-empty, 
which is what is required to link the routes.

\subsection{Local routes}
\label{sec-local}
Recall that each $T \in \TT$ is the distance from the origin to a nearby site on the diagonal, say $v_T$.   
We will prove
\begin{Proposition}
\label{P11}
There exists a random variable $D$, with all moments finite, such that
for each $T \in \TT$, each 
city (if any) $\xi$ in the associated scaled parallelogram $S_{v_T}$,  
and each city $\xi^*$ in the unit square centered at $(0,0)$, 
\[ \ell(\xi^*,\xi) \le D . \]
\end{Proposition}

\proof
Take one $T \in \TT$, say $T(a,\nwarrow)$.
Take a sector that is two sectors away away from $a$, in this case sector $c$, 
and consider (Figure 5) the four doubly-infinite oriented paths 
associated with $T(a,\nwarrow), T(a,\searrow), T(c,\nwarrow), T(c,\searrow)$.
These define a path circuit around the origin, of length 
at most $B_1 \Sigma$, where 
$\Sigma := T(a,\nwarrow) + T(a,\searrow) + T(c,\nwarrow) + T(c,\searrow)$ 
and $B_1$ (and $B_2$ below) are constants. 
As argued below Proposition \ref{P10},
we can construct an associated route in the network, of length at most 
$B_2 \Sigma$ 
for some $B_2$.

\setlength{\unitlength}{0.17in}
\begin{picture}(14,27)(-14,-14)
\multiput(5,-5)(-1,1){12}{\circle*{0.2}}
\multiput(-7,-11)(1,5){4}{\line(0,1){3}}
\multiput(-7,-8)(1,5){4}{\line(1,2){1}}

\multiput(1,-10)(1,5){3}{\line(0,1){3}}
\multiput(0,-12)(1,5){4}{\line(1,2){1}}

\multiput(6,3)(-5,-1){3}{\line(-3,0){3}}
\multiput(8,4)(-5,-1){4}{\line(-2,-1){2}}

\multiput(8,-5)(-5,-1){4}{\line(-3,0){3}}
\multiput(10,-4)(-5,-1){4}{\line(-2,-1){2}}

\put(0.3,-0.1){(0,0)}
\put(-5.5,4.0){$v^a_{{\mbox{\tiny $\nwarrow$ }}}$}
\put(-4,4){\circle*{0.35}}
\put(-2,2){\circle*{0.35}}
\put(-1.8,2.5){$v^c_{{\mbox{\tiny $\nwarrow$ }}}$}

\put(2.3,-2){$v^a_{{\mbox{\tiny $\searrow$ }}}$}
\put(2,-2){\circle*{0.35}}
\put(4.6,-5.9){$v^c_{{\mbox{\tiny $\searrow$ }}}$}
\put(5,-5){\circle*{0.35}}
\end{picture}

{\bf Figure 5.}
{\small  
Constructing a route encircling the origin.
}

\vspace{0.12in}
\noindent
Looking back at Figure 2, route construction relative to the square grid, consider a path in $\LL_L $ passing through $v$ but not an adjacent (in $\LL_L$) vertex $w$; the associated network route in this region 
cannot come closer to $w$ than  distance $L/\sqrt{2}$.
It follows that the route constructed above encircles the origin and does not come closer to the 
origin than $\delta L$, for some
constant  $\delta > 0$,.  By taking $L$ large we may assume the route 
encircles the unit square centered at the origin. 
The encircling route, being at length at most $B_2 \Sigma$, must stay within the disc of radius $B_2 \Sigma$. 
By connectivity of the network, from any city $\xi^*$ in the unit square 
there is a route which meets the encircling route and then 
continues to $v^a_\nwarrow$. 
This route stays within the disc of radius $B_2 \Sigma$. 
We can therefore bound its length by the r.v. 
$D(B_2 \Sigma)$ in Lemma \ref{L7} below; 
combining (\ref{Tkbound}) with the conclusion of Lemma \ref{L7}
establishes Proposition \ref{P11}.
\qed

At the end of the proof we used the following  straightforward  lemma.
\begin{Lemma}
\label{L7}
For the PP (with a city planted at the origin) and for $\sigma > 0$, define a r.v. $D(\sigma)$ to be the sum of Euclidean distances between all pairs of cities within the disc of radius $\sigma$.  Then for any r.v. $\Sigma$ with all moments finite, the r.v. $D(\Sigma)$ has  all moments finite.
\end{Lemma}

\subsection{Completing the proof}
\label{sec-completing}
Theorem \ref{T1} and Corollary \ref{C1} follow 
from Propositions \ref{P10} and  \ref{P11},  and (\ref{Tkbound}, \ref{dvdv}), as we now explain.

Consider cities $\xi^*$ and $\xi^{**}$ in the unit squares centered centered at the origin and at $v \in \Ints^2_L$.  
Proposition \ref{P11}, applied at the origin and at $v $, together with Proposition \ref{P10}, implies
\begin{equation}
 \ell(\xi^*,\xi^{**}) \leq B( 1 + |v| + T_1 + T_2)+ D_1 + D_2 
 \label{www}
 \end{equation}
where $T_1$ and $T_2$ satisfy (\ref{Tkbound}) and $D_1$ and $D_2$ have the distribution of $D$ in Proposition \ref{P11}.
So, writing $U$ for the unit square, 
\[ \sum_{\xi^* \in U} \sum_{\xi^{**} \in v + U} \ell(\xi^*,\xi^{**}) \leq \left(  B( 1 + |v| + T_1 + T_2)+ D_1 + D_2 \right)   \ N_1 N_2 \] 
where $N_1$ (resp. $N_2$) is the number of points of the PP in $U$ (resp. $v+U$).  
The random variables on the right side have all moments finite, uniformly in $v$, and so 
\[ \Ex \sum_{\xi^* \in U} \sum_{\xi^{**} \in v + U} \ell(\xi^*,\xi^{**}) = O(|v|) \mbox{ as } |v| \to \infty . \]
This is enough to establish  Theorem \ref{T1} 
(here $v$ is restricted to $\Ints^2_L$, but the bound holds for all sufficiently large $L$).
Corollary \ref{C1} is derived in the same way after first taking the $k$'th power in (\ref{www}).

\section{The relative neighborhood graph}
\label{sec-RNG}
 The definition of the  RNG (relative neighborhood graph) can be rephrased as follows.
For two points $v, w$ in the plane, define the {\em lune} $A_{v,w}$ as the intersection of the two discs of radii $d(v,w)$ centered at $v$ and at $w$.  
Then

 \begin{center}
 \includegraphics[width=75mm]{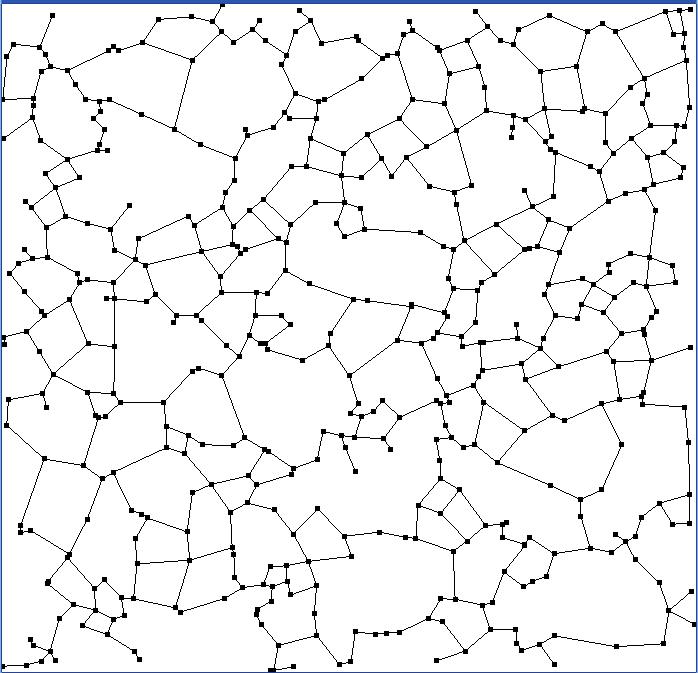}
 \end{center}
 
 {\bf Figure 6.} 
 {\small The relative neighborhood graph  on a realization of $500$ random cities.}
 
\vspace{0.1in}

\begin{eqnarray}
\mbox{there is a road between two cities $x,y$ if }\nonumber\\
\mbox{and only if  the lune  $A_{v,w}$ contains no other city.} \label{lune-d}
\end{eqnarray}
It is elementary \cite{jaromczyk} that the RNG contains the minimal spanning tree and is therefore connected (see section \ref{sec-MST} for the infinite case).

\begin{Proposition}
\label{P2}
The RNG on a PP satisfies the asymptotic essential connectedness property 
(\ref{AC}); and so by Theorem \ref{T1} it has the linearity property (\ref{lin-prop}).
\end{Proposition}
We'll discuss relations with previous results, after the proof.

\subsection{Proof of Proposition \ref{P2}}
Fix $L$. 
We start with a property for ``arbitrary" city positions -- this means ``deterministic, but in general position".  
\begin{Lemma}
\label{Ldd}
Consider an arbitrary finite configuration of cities in the $L \times L$ square.  Consider the network $\GG_L$ from (\ref{def-GG}); that is, the roads within the $L \times L$ square that are in the RNG for this configuration together with every locally finite configuration of cities outside the square. 
Suppose $\GG_L$ is not connected.  Write $\Delta(v)$ for the distance from $v$ to the boundary of the $L \times L$ square.  Then, for any initial city $v_0$, there exists a sequence of distinct cities 
$v_0, v_1,\ldots, v_m$, for some $m \geq 0$, such that \\
(i) $d(v_0,v_1) > d(v_1,v_2) > \ldots > d(v_{m-1},v_m)$\\
(ii) If $m = 0$ then every city  within distance $\Delta(v_0)$ from $v_0$ is in the same component as $v_0$\\
(iii) if $m \geq 1$ then $\Delta(v_m) < d(v_0,v_1)$.
\end{Lemma}
\proof
Write {\em components} for the connected components of $\GG_L$.
For each city $v$,  consider the nearest city ($w$, say) in some different component than $v$.
If the lune $A_{v,w}$  intersects the boundary of the $L \times L$ square, call $v$ a {\em terminal} city.
Otherwise, write $w = \theta(v)$.  Because $(v,w)$ is not an edge of the RNG, the lune $A_{v,w}$ must contain some other city, and 
(by definition of $w$ as {\em closest}) any such other city $z$ must be in the same component as $v$; note also $d(z,w) < d(v,w)$.
Thus if $w$ is not a terminal city then $\theta(w)$ satisfies 
$d(w,\theta(w)) \leq d(w,z) < d(v,w)$.
Finally, note that if $w = \theta(v)$ is a terminal city, then the lune $A_{w,u}$ intersects the boundary of the $L \times L$ square, 
where $u$ is the closest city to $w$ in a different component than $w$; since $d(w,u) \leq d(w,v)$ it must be that $w$ is within distance $d(v,w)$ from the boundary.

Combining the properties established above, we see that,
for any $v_0$, the sequence $v_1 = \theta(v_0), v_2 = \theta(v_1), \ldots, v_i = \theta(v_{i-1}), \ldots $ continued until reaching some terminal $v_m$, 
satisfies the assertions of the lemma.
\qed

We next give a lemma about Poisson points, not involving the RNG.
\begin{Lemma}
\label{L1}
Consider the points $\Xi$ of a PP in a $L \times L$ square. \\
(a) There is a constant $c^*$ such that, with probability $\to 1$  as $L \to \infty$, 
the length of the longest edge of the minimal spanning tree on the points $\Xi$  is at most
$c^* \log L$. \\
(b)  Fix $d_0$ and $n$.
The probability that there exists a sequence 
$\xi_0, \xi_1, \ldots , \xi_n$ of distinct points in $\Xi$ such that 
$d_0 \geq d(\xi_0,\xi_1) \geq d(\xi_1,\xi_2) \geq \ldots \geq d(\xi_{n-1},\xi_n)$ 
is at most 
$\frac{L^2 \pi^n d_0^{2n}}{n!}$.
\end{Lemma}
\proof Part (a) is a weaker version of the known precise asymptotics \cite{MR1442317} for longest MST edge length.  For (b), 
if we consider the PP on the whole plane, and require that $\xi_0$, but not the other $\xi_i$,  be in the $L \times L$ square, 
then the expected number of  sequences satisfying the condition in (b) equals
\[   L^2  \ \int \int  \ldots \int_{d_0 \geq r_1 \geq r_2 \ldots \geq r_n \geq 0}  
 2 \pi r_1 dr_1 \ 2 \pi r_2 dr_2 \ldots  2 \pi r_n dr_n ,
\]
so this is an upper bound on the desired probability.  
By symmetry of the $n!$ orderings, 
\[  \int \int  \ldots \int_{d_0 \geq r_1 \geq r_2 \ldots \geq r_n \geq 0}  
  r_1 dr_1 \  r_2 dr_2 \ldots  r_n dr_n  
  = \frac{1}{n!} \prod_{i=1}^n \int_0^{d_0} r_i dr_i 
  = \frac{1}{n!}  \ \left(\frac{d_0^2}{2}\right )^n \]
  and the result follows.
  \qed

We now start to combine the deterministic and random ingredients.
Write $\GG_L$ for the subnetwork (\ref{def-GG}) of the RNG on the points $\Xi$ of a PP in $[0,L]^2$.
\begin{Lemma}
\label{L2}
Fix $d_0$ and $n$.
The probability of the following event is at most 
$\frac{L^2 \pi^n d_0^{2n}}{n!}$.

There exist vertices $\xi^{\prime}$ and $\xi^{\prime \prime}$, in different components of $\GG_L$, with $\Delta(\xi^\prime) \geq d_0(n+1)$, such that $d(\xi^{\prime},\xi^{\prime \prime}) \leq d_0$.
\end{Lemma}
\proof  It is enough to show the event here implies the event in Lemma \ref{L1}(b).  We write (i,ii,iii) for the assertions of Lemma \ref{Ldd}.

Suppose such vertices $\xi^{\prime},\xi^{\prime \prime}$ exist.  Then $\GG_L$ is not connected, and we may consider the sequence $\xi^\prime = v_0, v_1, \ldots,v_m$ 
given by Lemma \ref{Ldd}.  By (ii) we have $m \ge 1$.
Because $d(v_0,v_1) \leq d(\xi^{\prime},\xi^{\prime \prime}) \leq d_0$, monotonicity in (i) implies $d(v_0,v_m) \leq d_0 m$. 
Because $\Delta(v_0) \le d(v_0.v_m) + \Delta(v_m)$, 
 we can apply (iii) to get 
$\Delta(v_0) \leq d_0 (m+1)$.   But by hypothesis $\Delta(v_0)  \geq d_0(n+1)$, so we have shown $m \geq n$. 
So the event in Lemma \ref{L1}(b) occurs.
\qed

\paragraph{Proof of Proposition \ref{P2}} 
Set
\[ n = n(L) = \lfloor L^{1/2} \rfloor, \  d_0 = d_0(L) = c^* \log L, \ M = M(L) = L - 2(1+ n(L))d_0 \]
and suppose $L$ is sufficiently large that $M > 1$.  As above, 
write $\GG_L$ for the subnetwork  (\ref{def-GG}) of the RNG on the points $\Xi$ of a PP in $[0,L]^2$. 
Consider the events

($A_L$): There exist vertices $\xi^{\prime}$ and $\xi^{\prime \prime}$, in different components of $\GG_L$, with $\Delta(\xi^\prime) \geq d_0(n+1)$, such that $d(\xi^{\prime},\xi^{\prime \prime}) \leq d_0$.

($B_L$): In the minimal spanning tree on the points of $\Xi$ within the $M \times M$ square concentric with $[0,L]^2$, the longest edge is longer than $d_0$.

Lemma \ref{L1}(a) implies $P(B_L) \to 0$ as $L \to \infty$.  Lemma \ref{L2} says 
\[ P(A_L) \leq \frac{L^2 \pi^n d_0^{2n}}{n!} \to 0 \mbox{ as } L \to \infty , \]
the convergence by a routine use of Stirling's formula and the definitions of $d_0(L)$ and $n(L)$.  
Now we assert 
\begin{eqnarray}
\mbox{if neither $A_L$ not $B_L$, then all the cities in the concentric} \nonumber\\
\mbox{ $M \times M$ square are in the same component of $\GG_L$.} 
\label{assert}
\end{eqnarray}
For if the implication were false, then some edge $(\xi^\prime, \xi^{\prime \prime})$ of the MST 
(on cities in the $M \times M$ square) links cities in different components of $\GG_L$, and because $B_L$ fails we have
$d(\xi^\prime, \xi^{\prime \prime}) \leq d_0$.
But $\Delta(\xi^\prime) \ge   (L-M)/2 = d_0(n+1)$ by definition of $M$,
 so  $A_L$ holds.

But assertion (\ref{assert}) implies 
\[ N^0_L  \leq N_{L \setminus M} + N_L 1_{A_L \cup B_L} \]
where $N^0_L$ is the number of
cities inside $[0,L]^2$ that are not in the largest component of $\GG_L$,
 $N_L$ is the total number of
cities inside $[0,L]^2$, and  $N_{L \setminus M}$ is the number of cities  inside $[0,L]^2$ but not inside  the concentric $M \times M$ square.
Taking expectation, 
\[ L^{-2} \Ex N^0_L \leq 1 - \sfrac{M^2}{L^2} + \Ex (L^{-2}N_L \  1_{A_L \cup B_L} )\] 
and the right side $\to 0$ by definition of $M(L)$, uniform integrability of $(L^{-2}N_L)$ and the fact 
$P(A_L \cup B_L) \to 0$. 
\qed

\subsection{Proximity graphs}
\label{sec-discuss}
There is a general class of {\em proximity graphs} \cite{jaromczyk} which as defined at (\ref{lune-d}) but with the lune $A_{v,w}$ replaced by some specified 
subset of the lune; such  networks on a PP {\em a priori} have more roads than the RNG and therefore can only have shorter network distances.  
So the linearity property remains true for this class of networks.

For any finite spatial network one can define its {\em stretch} as the maximum, over city-pairs $(x,y)$, of the ratio 
$\ell(x,y)/d(x,y)$.
Now consider a network, in the sense of a deterministic rule for assigning roads to any  finite configuration of cities.  
Such a network is called \cite{MR2289615} a {\em $t$-spanner} if the stretch is always bounded by the constant $t$.  
Clearly the property of being a {\em spanner} (that is, a $t$-spanner for some $t$) is stronger than our  linearity property (\ref{lin-prop}).
It is a remarkable result 
\cite{MR1134449}
that the Delaunay triangulation is a $t$-spanner for 
$t = \frac{2\pi}{3 \cos \pi/6} \approx 2.42$.  
Consequently, while our Theorem \ref{T1} can easily be applied to the Delaunay triangulation over the PP, it does not yield any new result.  
On the other hand the RNG and related proximity graphs are known {\em not} to be  spanners, and (when constructed over $n$ random points in a square) the expectation of stretch 
increases slowly to infinity as $n \to \infty$.  See \cite{MR2257270}.  
So our results do say something new about this model.

\subsection{Modifications of disconnected networks}
Our results suggest the following general program,
which we have not pursued.
Starting with (for instance) the supercritical geometric random graph, one might guess that any reasonable {\em ad hoc} scheme for adding edges to create a connected network would yield a network with the linearity property, and one could seek to verify this for any given scheme by using Theorem \ref{T1}.
Similarly, starting from a more general disconnected network, one can always create a connected network by adding all edges $(v,w)$ of the RNG such that $v$ and $w$ are in different components of the original network; given any particular original network one could seek to modify the proof of Proposition \ref{P2} to show that 
 the asymptotic essential connectedness property 
(\ref{AC}) holds and thereby establish  the linearity property.

\subsection{Remarks about the minimum spanning tree}
\label{sec-MST}

On a general infinite (but locally finite) configuration  of points, the analog of the MST is the {\em minimum spanning forest}.  It is true  \cite{alexander95}, but not obvious, that when applied to the PP one gets a single tree.  
This implies in particular that the RNG over the PP is connected.  

It is widely believed in the statistical physics literature (see e.g.  \cite{read-2005-72,wieland} for references) that for the MST over the PP, quantities such as $\Ex T_z$ should grow as some power $|z|^\gamma$ where $\gamma \approx 1.23$.  
But it is not clear to us what has been rigorously proved.  
Intuitively, \underline{no} tree network
can have the linearity property (\ref{lin-prop}), and proving this would be a minor research project.


\end{document}